\documentclass{article}
\usepackage[utf8]{inputenc}
\usepackage{amsmath}
\usepackage{amsfonts}
\usepackage{amssymb}
\usepackage{mathtools}
\usepackage{amsthm}

\title{Locally Maximizing Metric of Width on Manifolds with Boundary}
\author{Yucheng Tu}
\date{}

\begin{document}

\maketitle
\begin{abstract}
    In this paper we use min-max theory to study the existence free boundary minimal hypersurfaces (FBMHs) in compact manifolds with boundary $(M^{n+1}, \partial M, g)$, where $2\leq n\leq 6$. Under the assumption that $g$ is a local maximizer of the width of $M$ in its comformal class, and all embedded FBMHs in $M$ are properly embedded, we show the existence of a sequence of properly embedded equidistributed FBMHs. This work extends the result of Ambrozio-Montezuma [2]. 
\end{abstract}
\vspace{0.3in}

\section{Introduction}
In a recent work of Ambrozio and Montezuma [2], the equidistribution phenomenon of minimal $S^2$ in $S^3$ is studied. With the assumption that the metric $g_0$ on $S^3$ is a local maximizer(in its conformal class) of the Simon-Smith width functional $W(S^3,g)$, the authors proved the existence of equi-distributed minimal 2-spheres in measure theoretic sense. In this paper we follow their main ideas and extend the results to embedded free boundary minimal hypersurface in a ball of dimension $3\leq n+1\leq 7$. We shall prove the following result:\\

\noindent\textbf{Theorem 1.1} Given metric $g$ on $(M^{n+1},\partial M)$, $2\leq n\leq 6$, if $g$ maximizes the normalized width $W(M,g)$ in the conformal class of $g$, and suppose that all free boundary minimal hypersurfaces embedded in $M$ is properly embedded, then there exist a sequence $\{\Sigma_i^n\}$ of free boundary minimal hypersurfaces with index zero or one and area no greater than $W(M,g)$ for which the following holds for all $f\in C(M)$:
$$\lim_{k\rightarrow\infty}\frac{1}{\sum_{i=1}^{k}\text{area}(\Sigma_i,g)}\sum_{i=1}^k\int_{\Sigma_i}fdA_g=\frac{1}{\text{vol}(M,g)}\int_{M}fdV_g.$$

Furthermore, if we assume that $(M,\partial M,g)$ contains no stable free boundary minimal hypersurface with area greater than its width $W(M,g)$, then we can drop the proper embeddedness assumption, and choose $\{\Sigma_i\}$ so that each of them has index $1$ and area equal to $W(M,g)$:\\

\noindent\textbf{Theorem 1.2} Given metric $g$ on $(M^{n+1}, \partial M)$, $2\leq n\leq 6$, if $g$ maximizes the normalized width in the conformal class of $g$ and there exists no stable free boundary minimal hypersurface of area less than $W(M,g)$, then there exist a sequence $\{\Sigma_i^n\}$ of free boundary minimal hypersurfaces with index one and area equal to $W(M,g)$ for which the following holds for all $f\in C(M)$:
$$\lim_{k\rightarrow\infty}\frac{1}{kW(M,g)}\sum_{i=1}^k\int_{\Sigma_i}fdA_g=\frac{1}{\text{vol}(M,g)}\int_{M}fdV_{g}.$$

The main difference between our theorem and Proposition 1.4.1 of [2] is that in free boundary case, we can not rule out the case when $\{\Sigma_i\}$ is not properly embedded ($\Sigma\cap \partial M\neq\emptyset$), due to the lack of convexity of $\partial M$.  Readers can see [6] for a possible example of non-properly embedded free boundary minimal hypersurface in an Euclidean domain.\\

\noindent\textbf{Acknowledgment.} The author would like to thank Prof. Xin Zhou for introducing and explaining the topic to him, and very helpful discussions.

\section{Preliminaries}

In the following let $2\leq n\leq 6$, and $(M^{n+1},\partial M, g)$ be a Riemannian manifold with smooth boundary $\partial M$ and metric $g$. The notions of sweepout and width are crucial in the min-max theory of minimal hypersurfaces. In [6] a min-max theory of free boundary minimal hypersurfaces(FBMH) were developed, which is of great use here in our context. First we give a introduction to FBMH. 

\subsection{Free Boundary Minimal Hypersurfaces and Morse Index}
Let $(M,\partial M, g)$ be as above. A free boundary minimal hypersurface $\Sigma$ in $(M,g)$ is a $n$-dimensional submanifold of $M$ with vanishing mean curvature($H=0$) and boundary $\partial\Sigma$ orthogonal to $\partial M$. We can also use the first variation of area of $\Sigma$ to characterize this property: given a smooth perturbation of $M$ defined by $\phi:M\times(-\epsilon,\epsilon)\rightarrow M$ with $\phi(\cdot,0)=\text{id}_M$ and $\phi(\partial M,\cdot)\subset\partial M$, we have the following first variation formula:
$$\frac{\partial}{\partial s}\text{area}[\phi(\Sigma,s)]\bigg|_{s=0}=\int_\Sigma -H\vec{n}\cdot \frac{\partial\phi}{\partial t}\Big|_{s=0}dA+\int_{\partial\Sigma}\phi\eta\cdot \vec{n}ds$$
where $\vec{n}$ is the unit normal of $\Sigma$ and $\eta$ is the outward conormal along $\partial\Sigma$. Therefore $\Sigma$ is a critical point if and only if $H=0$ on $\Sigma$ and $\eta\perp \vec{n}$ on $\partial\Sigma$, as in the definition of FBMS. For variation in normal direction as $\frac{\partial\phi}{\partial t}\Big|_{s=0}=f\cdot\vec{n}$, we have the second variation of area:
$$\frac{\partial^2}{\partial s^2}\text{area}[\phi(\Sigma,s)]\bigg|_{s=0}=\int_{\Sigma}\big(|\nabla f|^2-\text{Ric}_M(\vec{n},\vec{n})f^2-|A|^2|f|^2\big)d\mu-\int_{\partial\Sigma}h^{\partial M}f^2ds$$
After an integration by part, the right hand side of the second variation formula defines a quadratic form on $C^\infty(\Sigma)$:
$$I(f,g)=\int_\Sigma(-f\Delta g-\text{Ric}_M(\vec{n},\vec{n})fg-|A|^2fg)d\mu+\int_{\partial\Sigma}\big(f\frac{\partial g}{\partial\vec{n}}-h^{\partial M}fg\big)ds$$
and we define the index of $\Sigma$ to be the number of negative eigenvalues of $I$. $\Sigma$ is called a stable FBMS if its index is 0, i.e. there is no variation that reduce the area of $\Sigma$ to the second order.
\subsection{Min-max Construction}
Given a manifold with boundary $(M^{n+1},\partial M,g)$, let $\mathcal{Z}_n(M,\partial M,\mathbb{Z})$ be the space of integer rectifiable $n$-currents $T$ in $M$ with coefficients in $\mathbb{Z}$, such that $\partial T\in\partial M$, modulo the following equivalence relation:
$$T\sim S\text{ iff } T-S\in\mathcal{R}_n(\partial M, \mathbb{Z})$$
where $\mathcal{R}_n(\partial M, \mathbb{Z})$ is the space of $n$-rectifiable integral currents in a sufficiently high dimensional Euclidean space $\mathbb{R}^L$, supported on $\partial M$. (We can regard the $M$ as embedded isometrically in $\mathbb{R}^L$.) We endow $\mathcal{Z}_n(M,\partial M,\mathbb{Z})$ with the flat topology $\mathcal{F}$. Let us define the notion of 1-sweepout and 1-width.\\

\noindent\textbf{Definition 2.1}(cf [5]) Let $(M,\partial M)$ be defined as above. A 1-sweepout of $M$ is a one parameter family of maps $\Phi:[-1,1]\to \mathcal{Z}_n(M,\partial M,g)$ satisfying the following conditions:
\begin{itemize}
    \item[(1)] $\Phi$ is continuous in flat topology;
    \item[(2)] $\sup_{x\in I}\mathbf{M}(\Phi(x))<+\infty$;
    \item[(3)] there is no mass concentration on $\Phi$;
    \item[(4)] $F(\Pi_\Phi)$ represents a non-zero element in $H_{n+1}(M,\partial M)$.
\end{itemize}

\noindent \textbf{Definition 2.2} We define the width of a manifold with metric $g$ as $$W(M,\partial M,g)=\inf_{\Phi\in\Bar{\Lambda}}\Big(\max_{t\in[-1,1]}\mathbf{M}(\Phi(t),g)\Big)$$
where $\Phi$ is a sweepout of $(M,\partial M,g)$. The normalised width is defined by $$W_n(M,\partial M,g)=\frac{W(M,\partial M,g)}{\text{Vol}(M,g)^{\frac{n}{n+1}}}.$$
Let us note that by a similar argument as in [7], under a smooth variation of metric $g(t)$ with respect to the original metric, $W(M,g(t))$ is a Lipshitz function of $t$.

\section{Proof of the Main Theorems}

In this section we prove Theorem 1.1 using a perturbation method originally due to Marques-Neves-Song[7], and prove Theorem 1.2 by a calculation of derivative of width inspired by Fraser-Schoen's work[3] on Steklov eigenvalues.

\subsection{Proof of Theorem 1.1}

In view of the abstract theorem 4.2, we can reduce the equi-distribution property to proving the following lemma:\\

\noindent\textbf{Lemma 3.1.1} Let $g$ be a Riemannian metric on $M$ that maximizes the normalized width in its conformal class. For every continuous function $f$ satisfying 
$$\int_{M}fdV_g<0,$$
there exists some integers $n_1,\cdots, n_N$, and disjoint embedded free boundary minimal hypersurfaces $\Sigma_1,\cdots,\Sigma_N$ in $(M,g)$ such that 
$$W(M,g)=\sum_{i=1}^{N}n_j\text{area}(\Sigma_i,g),\quad \sum_{i=1}^N\text{Ind}_g(\Sigma_i)\leq 1$$
and
$$\sum_{i=1}^N n_i\int_{\Sigma_i}fdA_g\leq 0.$$

In order to associate the function $f$ with the derivative of width under a conformal change of metric, we need to perturb the conformal family of the original metric to a new family so that the width is differentiable. The following technical lemma is crucial:\\

\noindent\textbf{Lemma 3.1.2} Let $q\geq 4$ be an integer, and $g:[0,1]\to\Gamma_q$ be a smooth embedding. Then there exist smooth embeddings $h:[0,1]\to\Gamma_q$ which are arbitrarily close to $g$ in the smooth topology, and $J\subset[0,1]$ with full Lebesgue measure such that
\begin{itemize}
    \item[(1)] The function $W(M,h(t))$ is differentiable at every $\tau\in J$; and
    \item[(2)] For each $\tau\in J$, there exist a collection of integers $\{n_1,\cdots, n_N\}$ and a finite collection $\{\Sigma_1,\cdots,\Sigma_N\}$ of disjoint free boundary embedded minimal hypersurfaces of class $C^q$ in $(M,h(\tau))$ such that
    $$W(M,h(\tau))=\sum_{k=1}^Nn_k\cdot\text{area}(\Sigma_k,h(\tau)),\quad \sum_{k=1}^N\text{ind}_{h(\tau)}(\Sigma_k)\leq 1,$$
    $$\text{and}\quad \frac{d}{dt}\Big|_{t=\tau}W(M,h(t))=\frac{1}{2}\sum_{k=1}^Nn_k\int_{\Sigma_k}\text{Tr}_{(\Sigma_k,h(\tau))}(\partial_th(\tau))dA_{h(\tau)}.$$
\end{itemize}
\begin{proof}
(cf [2]) First, due to the density of bumpy metric on $M$ and Rademacher's theorem, we can perturb the smooth family $g:[0,1]\to\Gamma_q$ to $h:[0,1]\to\Gamma_q$ which is arbitrarily close to $g$ in smooth topology, and a set $J\subset[0,1]$ of full measure such that $h(\tau)$ is a bumpy metric and $W(M,h(t))$ is differentiable at $\tau$, for all $\tau\in J$.

For all $\tau\in J$, fix a sequence $t_i\to \tau$, we have 
$$\frac{d}{dt}W(M,h(t))\Big|_{t=\tau}=\lim_{i\to\infty}\frac{W(M,h(t_i))-W(M,h(\tau))}{t_i-\tau}.$$

By Proposition 7.3 of [4], we can find a finite disjoint collection of FBMHs $\{\Sigma_1(t_i),\cdots,\Sigma_{i_k}(t_i)\}$ and integers $\{N_1,\cdots, N_{i_k}\}$ such that
$$W(M,h(t_i))=\sum_{j=1}^k N_j\text{area}(\Sigma_{i_j}(t_i))\qquad \sum_{j=1}^k N_j\cdot\text{Ind}(\Sigma_{i_j}(t_i))\leq 1$$
Now as $t_i\to \tau$, since $h$ is a smooth family we have $\text{area}(\Sigma_{i_j}(t_i))$ uniformly bounded below and above by $W(M,h(\tau))$ as $t_i$ is sufficiently close to $\tau$. Therefore by the compactness theorem we can extract a subsequence $t_{i_j}$ so that $\Sigma_{i_{j_k}}$ converges in the varifold sense to $\Sigma_k$, since the metric $h(\tau)$ is bumpy, there is no multiplicity issue in the convergence, so we can conclude that the convergence is graphical and smooth. Therefore standard calculation shows
$$\lim_{i\to\infty}\frac{\text{area}(\Sigma_{i_{j_k}},h(t_{i_j}))-\text{area}(\Sigma_j,h(\tau))}{t_{i_j}-\tau}=\frac{1}{2}\int_{\Sigma_k}\text{Tr}_{(\Sigma_k,h(\tau))}(\partial_th(\tau))dA_{h(\tau)}$$
and hence we have the derivative of width formula.
\end{proof}

\noindent Now we can finish the proof of Theorem 1.1 by showing Lemma 3.1.1. For a continuous function $f$ with $\int_{M}fdV_g<0$, we can define a conformal change of metric:
$$g(t)=(1+\frac{n+1}{n}tf)^\frac{n}{n+1}g\quad\text{for}\quad 0\leq t\leq T.$$
We have $\partial_tg(t)\big|_{t=0}=fg$, hence for small $T>0$ we have $\text{Vol}(M,g(t))$ less than the the volume under the original metric. Since $g$ maximizes the normalised width, we have
$$\frac{W(M,g(t))}{\text{Vol}(M,g(t))^{\frac{n}{n+1}}}\leq\frac{W(M,g(0))}{\text{Vol}(M,g(0))^{\frac{n}{n+1}}}\quad\text{for}\quad 0\leq t\leq T.$$
Hence 
$$W(M,g(t))\leq W(M,g(0)) \Big(\frac{\text{Vol}(M,g(t))}{\text{Vol}(M,g(0))}\Big)^{\frac{n}{n+1}}< W(M,g(0))\quad\text{for}\quad 0\leq t\leq T.$$
Fix $q\geq 4$. Now for each $i\in\mathbb{N}$ with $1/i<T$, we can find a perturbation $h_i:[0,1/i]\to\Gamma_q$ and $J_i\subset[0,1/i]$ with full Lebesgue measure such that
$$W(M,h_i(1/i))< W(M,h_i(0))$$ 
and so there is $\tau_i\in J_i$ such that
$$\frac{d}{dt}W(M,h_i(t))\Big|_{t=\tau_i}\leq 0$$
due to the first fundamental theorem of calculus. Hence by the previous lemma there are FBMHs $\Sigma_{i_j}$, $j=1,2,\cdots, n_i$ and a set of integers $\{n_{i_1},\cdots,n_{i_N}\}$ such that 
$$W(M,h_i(\tau_i))=\sum_{k=1}^Nn_{i_k}\cdot\text{area}(\Sigma_{i_k},h_i(\tau_i)),\quad \sum_{k=1}^N\text{ind}_{h_i(\tau_i)}(\Sigma_{i_k})\leq 1,$$
$$\text{and}\quad \frac{d}{dt}\Big|_{t=\tau_i}W(M,h_i(t))=\frac{1}{2}\sum_{k=1}^Nn_{i_k}\int_{\Sigma_{i_k}}\text{Tr}_{(\Sigma_{i_k},h_i(\tau_i))}(\partial_th_i(\tau_i))dA_{h_i(\tau_i)}\leq 0.$$
We can relabel these $\Sigma_{i_k}$ such that except for $\Sigma_{i_1}$, others have index $0$. Now we can use the Compactness Theorem A.6 to conclude that, by picking a subsequence $\tau_{i_j}\to 0$, the FBMHs subconverges smoothly and graphically to $\{\Sigma_1,\cdots,\Sigma_N\}$ with multiplicity $1$, except for $\Sigma_1$, where the multiplicity can be $2$ if $\Sigma_1$ is stable. Therefore we can pass the limit of the formula above and show
$$\frac{1}{2}\sum_{k=1}n_k\int_{\Sigma_k}fdA_g=\frac{1}{2}\sum_{k=1}n_k\int_{\Sigma_k}\text{Tr}_{\Sigma_k}(\partial_tg(0))dA_g\leq 0$$
Hence this finish the proof when $f$ is a smooth function on $(M,g(0))$. When $f$ is a continuous function we can use smooth functions to approximate $f$ uniformly and use similar arguments in the sequence picking process. Hence once Lemma 3.1.1 is proved, then the implication i) to iv) in \textbf{Theorem 4.2} applies if we let $Y$ be the Radon measure induced from FBMHs in $M$, and $\mu_0$ be the Hausdorff measure on $(M,g)$. 

\subsection{Proof of Theorem 1.2}
Now we prove Theorem 1.2. First we need a result that guarantees the existence of optimal sweepout in Lemma 1.2.2, and then we can compute the derivative of width under a general smooth family of metrics.\\

\noindent\textbf{Lemma 3.2.1} ([5] Prop. 5.4). Let $(\Sigma,\partial\Sigma)\subset (M,\partial M)$ be an orientable, almost properly embedded, free boundary minimal hypersurface with $\text{Area}(\Sigma)$ less than the area of the stable free boundary minimal hypersurface in $M$. Then there is a sweepout 
$$\Psi: [-1,1]\rightarrow Z_n(M,\partial M),$$
such that:\\
(1) $\Psi(0)=\Sigma$;\\
(2) $F(\Psi)=M$;\\
(3) $\mathbf{M}(\Psi(t))<\text{Area}(\Sigma)$ for $t\neq 0$.\\

\noindent\textbf{Lemma 3.2.2} Given $\{g(t)\}_{t\in (a,b)}$ as a one parameter family of metrics on $M$ varying smoothly, if $t_0\in (a,b)$ is a point where $W(t):=W(M,g(t))$ is differentiable, then there is an almost properly embedded free boundary minimal hypersurface $\Sigma$ in $(M,g(t_0))$ such that
$$\text{area}(\Sigma,g(t_0))=W(t_0)\quad\text{and}\quad\frac{d}{dt}W(M^n,g(t))\Big|_{0}=\frac{1}{2}\int_{\Sigma}\text{Tr}_{\Sigma}(\frac{\partial}{\partial t}g(t)\Big|_{t=t_0})dA_{g(t_0)}.$$
\begin{proof}
By Lemma 1.2.1, there exist an optimal sweepout $\{\Sigma_s\}_{s\in[-1,1]}$ such that $\text{area}(\Sigma_{0})=W(M,g(t_0))$ and for all $s\neq 0$, $\text{area}(\Sigma_s)<\text{area}(\Sigma_0)$. Consider a smooth function $F:(a,b)\times [-1,1]\rightarrow \mathbb{R}$ defined as $F(t,s)=\text{area}(\Sigma_s, g(t))$, then we have $F_s(t_0,0)=0$ and $F_{ss}(t_0,0)< 0$. Now let us show that there exists $\epsilon>0$ such that there is a differentiable function $s=s(t)$ for $t\in (t_0-\epsilon,t_0+\epsilon)$, such that
$$F(t,s(t))=\max_{s\in[-1,1]}F(t,s).$$
Since $F_{ss}(t_0,0)<0$, the implicit function theorem guarentees that $F_s(t,s)=0$ defines a smooth function $s=s(t)$ on $(t_0-\epsilon, t_0+\epsilon)$. Now there is a neighborhood of $(t_0,0)$ such that $F_{ss}<0$, and therefore $F(t,s(t))$ is a local maximum for each fixed $t\in(t_0-\epsilon',t_0+\epsilon')$. Due to the construction of sweepout(property 3) and possibly making $\epsilon'$ even smaller we can make sure $F(t,s(t))$ is a strict maximum. Hence the claim is proved. Now we define a function $h(t)=F(t,s(t))-W(t)$ over a neighborhood of $t_0$. We have that $h(t)\geq 0$ due to the definition of width, and $h(t_0)=0$ is the local minimum. Since $W(t)$ is differentiable at $t_0$, $h$ is also differentiable and $h'(t_0)=0$. Hence we have
$$W'(t_0)=\frac{\partial}{\partial t}F(t,s(t))\big|_{t=t_0}=F_s(t_0,0)s'(t_0)+F_t(t_0,0)=\frac{1}{2}\int_{\Sigma}\text{Tr}_{\Sigma}(\frac{\partial g}{\partial t}(t_0))dA_{g(t_0)}$$
\end{proof}

Similar to the proof of Theorem 1.1, we can define a conformal change of the metric $g$, this time with a volume preserving factor. More precisely, for a smooth function $f$ with $\int_{\mathbb{B}}fdV_g=0$, we fix a small $T>0$ and let
$$g(t)=\frac{\text{Vol}(M,g)^{\frac{n}{n+1}}(1+ft)}{\text{Vol}(M,(1+ft)g)^{\frac{n}{n+1}}}g\quad\text{for all }t\in [0,T).$$
It is straightforward to show that $\text{Vol}(M,g(t))=\text{Vol}(M,g(0))$ for all $t\in[0,T)$, and that $\partial_tg(0)=fg$.\\

\noindent\textbf{Lemma 3.2.3} Let $g(t)$, $t\in [0,\epsilon)$ be a smooth family of Riemannian metrics on $M$ that contains no stable free boundary minimal surface with area greater than $W(M,g)$. If 
$$W(M,g(0))\geq W(M,g(t))$$
then there exists a free boundary minimal disk $\Sigma$ such that
$$\text{area}(\Sigma,g(0))=W(M,g(0))\quad\text{and}\quad\int_{\Sigma}\text{Tr}_\Sigma(\partial_tg(0))dA_{g(t_0)}\leq 0.$$
\begin{proof}
Take an $\epsilon>0$. By Rademacher's Theorem, $W$ is differentiable at almost all $t\in[0,\epsilon)$. Since $W$ assumes local maximum at $0$, There exists a sequence $t_n\in[0,\epsilon)$ converging to $t_0$ such that $W'(t_n)\leq 0$ for all $n$. Hence by the previous lemma we can find an embedded free boundary minimal disk $\Sigma_n$ in $(M,g(t_n))$ with $\text{area}(\Sigma_n,g(t_n))=W(t_n)$ and $\int_{\Sigma_n}\text{Tr}_{\Sigma_n}(\partial_tg(t_n))dA_{g(t_n)}\leq 0$. Now by the compactness theorem we see that $\Sigma_n$ subconverges to a embedded free boundary minimal disk $\Sigma$. By the smooth convergence we have $\text{area}(\Sigma,g(0))=W(0)$ and $\int_\Sigma\text{Tr}_\Sigma(\partial_tg(0))dA_{g(0)}\leq 0$.  
\end{proof}

Combining the Lemma 1.2.3 and the previously defined conformal change of metric, we can show the following statement:\\

\noindent\textbf{Proposition 3.2.4} Let $f$ be a continuous function on $(M,g)$ with zero average, and if $(M,g)$ contains no stable free boundary minimal surface with area greater than $W(M,g)$, we can find a almost properly embedded FBMH $\Sigma$ in $(M,g)$ such that $\int_\Sigma fdA_g\leq 0$.
\begin{proof}
This statement follows when we approximate the function $f$ uniformly by smooth functions, and use the previous conformal change of metric.
\end{proof}

Then as in the proof of Theorem 1, the implication ii) to iv) in \textbf{Theorem 4.2} will confirm the existence of equidistibuted FBMHs in $M$, and as \textbf{Lemma 3.2.3} shows, each $\Sigma_i$ has area equal to $W(M,g(0))$.

\section{Compactness Theorem and Equidistribution Theorem}

In this section we prove a compactness theorem of FBMH for varying background metric and the abstract theorem on the existence of equi-distributed sequence of measures.\\ 

\noindent\textbf{Theorem 4.1.} Let $2\leq n\leq 6$ and $N^{n+1}$ be a compact manifold with boundary and $\{g_k\}_{k\in\mathbb{N}}$ a family of Riemannian metrics on $N$ converging smoothly to some limit $g$. If $\{M_k^n\}\subset N$ is a sequence of connected and properly embedded free boundary minimal hypersurface in $(N,g_k)$ with 
$$H^n(M_k)\leq \Lambda<\infty\quad\text{and}\quad \text{index}_k(M_k)\leq I,$$
for some fixed constants $\Lambda\in\mathbb{R}$, $I\in\mathbb{N}$, both independent of $k$. Then up to subsequence, there exists a connected and free boundary embedded minimal hypersurface $M\subset (N,g)$ where $M_k\rightarrow M$ in the varifold sense with
$$H^n(M)\leq \Lambda<\infty\quad \text{index}_k(M_k)\leq I$$
we have that the convergence is smooth and graphical for all $x\in M-Y$ where $Y=\{y_i\}_{i=1}^K\subset M$ is a finite set with $K\leq I$ and the following dichotomy holds:
\begin{itemize}
    \item if the number of leaves in the convergence  is one then $Y=\Phi$, i.e.
    the convergence is smooth an graphical everywhere
    \item if the number of sheets is $\geq 2$\\
    -if $N$ has $\text{Ric}_N>0$ then $M$ cannot be one-sided\\
    -if $M$ is two-sided the $M$ is stable.
\end{itemize}

\begin{proof}
We know by Allard's compactness theorem that there is an $M$ such that after passing to a subsequence, $M_k\rightarrow M$ in $\mathbf{IV}_n(N)$. Let $Y\subset M$ be the singular set of $M$. First we show that $|Y|\leq I$. Suppose on the contrary that $Y$ contains at least $I+1$ points $y_1,\cdots, y_{I+1}$. Then we can find $\{\epsilon_i\}_{1}^{I+1}$ such that $B(y_i,\epsilon_i)\cap B(y_j,\epsilon_j)=\emptyset$, and that $\sup_{k}\sup_{M_k\cap B(y_i,\epsilon_i)}|A|^2=\infty$, for all $i=1,\cdots,I+1$. Since $g_k$ converges to $g$ smoothly, the sectional curvature of $(N^{n+1},g_k)$ are uniformly bounded. Hence curvature estimate of [4] applies to this varying metric case, that is, in $\Sigma_k\cap B_r(p)$ the second fundamental form of $\Sigma_k$ are bounded by a uniform constant $C$ that depends only on $N$. Hence
we infer that for sufficiently large $k$, $M_k\cap B(y_i,\epsilon_i)$ is not stable for all $i=1,\cdots,I+1$. This implies that $\text{index}_k(M_k)\geq I+1$ which contradicts with the assumption.\\

To Show that $\text{Index}(M)\leq I$, we suppose that there are $u_1,u_2,\cdots,u_{I+1}\in C^\infty(M)$ that are $L^2$-orthogonal such that $I(u_i,u_i)<0$ for $i=1,2,\cdots,I+1$. Then we extend $u_i$ to $\Tilde{u}_i\in C^1(M)$ and let $u_i^{k}=\Tilde{u}_i\big|_{M_k}$. Since $M_k\rightarrow M$ as varifold, we have for sufficiently large $k$, $I_k(u_i^k,u_i^k)<0$ for $i=1,2,\cdots,I+1$. Since $\text{Index}(M_k)\leq I$, $\{u_i^k\}_{i=1}^{I+1}$ must be linearly dependent. By taking a subsequence and relabeling if necessary, we can find $\{\lambda_i\}_{i=1}^I\subset \mathbb{R}$ and $\lambda_i$'s not all zero such that $u_{I+1}^k=\sum_{i=1}^n\lambda_iu_i^k$. By varifold convergence we have $\langle u_i^k,u_j^k\rangle\rightarrow\langle u_i,u_j\rangle=\delta_{ij}$ for $i,j=1,2,\cdots,n+1$. Therefore by the varifold convergence,

$$0=\langle u_{n+1},u_i\rangle_M=\lim_{k\rightarrow\infty}\langle u^k_{n+1},u^k_i\rangle_{M_k}=\lim_{k\rightarrow\infty}\lambda_i$$
This implies that $u_{n+1}=0$ which contradicts $I(u_{n+1},u_{n+1})<0$.\\

Now if the multiplicity of convergence is $1$, then the convergence is smooth everywhere by the regularity theorem of [8]. Hence the theorem is proved.
\end{proof}

Here we include an abstract theorem used in the proof of Theorem 1.1 and 1.2., for the proof see Theorem B.1 of [2].\\

\noindent\textbf{Theorem 4.2.}(cf [2]) Let $Y$ be a non-empty weak-* compact subset of $M(X)$. The following assertions about a measure $\mu_0$ in $M(X)$ are equivalent to each other:
\begin{itemize}
    \item[i)]  For every function $f\in C^0(X)$ such that $\int_Xfd\mu_0<0$, there exists $\mu\in Y$ such that $\int_Xfd\mu\leq 0$.
    \item[ii)] For every function $f\in C^0(X)$ such that $\int_Xfd\mu_0=0$, there exists $\mu\in Y$ such that $\int_Xfd\mu\leq 0$.
    \item[iii)] $\mu_0$ belongs to the weak-* closure of the convex hull of the positive cone over $Y$.
    \item[iv)] There exists a sequence $\{\mu_k\}$ in $Y$ such that
    $$\lim_{k\rightarrow\infty}\frac{1}{k}\sum_{i=1}^{k}\frac{1}{\mu_i(X)}\int_{X}fd\mu_i=\frac{1}{\mu_0(X)}\int_Xfd\mu_0\quad\text{for all }f\in C^0(X).$$
\end{itemize}

\end{document}